\documentclass{amsart} \usepackage{hyperref}
\usepackage{enumerate} 
\usepackage{upref}
\usepackage{verbatim} 
\usepackage{color}
\usepackage{graphicx}
\usepackage{todonotes}
\newcommand{\sign}{\mathop{\rm sign}}
\newcommand*{\mailto}[1]{\href{mailto:#1}{\nolinkurl{#1}}}

\usepackage[latin1]{inputenc}
\usepackage{amssymb, amsmath, amsfonts, amsthm}
\usepackage[all]{xy}

\usepackage{marginnote}

\DeclareMathOperator{\id}{Id}

\DeclareMathOperator{\supp}{supp}
\newcommand{\dott}{\, \cdot\,}

\newcommand{\D}{\ensuremath{\mathcal{D}}}

\newcommand{\F}{\ensuremath{\mathcal{F}}}

\newcommand{\abs}[1]{\left\vert#1\right\vert}

\newcommand{\Real}{\mathbb R}

\newcommand{\norm}[1]{\left\Vert#1\right\Vert}

\newcommand{\muac}{\mu_{\text{\rm ac}}}

\newtheorem{theorem}{Theorem}[section]

\newtheorem{lemma}[theorem]{Lemma}
\newtheorem{definition}[theorem]{Definition}

\newtheorem{remark}[theorem]{Remark}
\numberwithin{equation}{section}

\allowdisplaybreaks

\begin{document}

\title[Wave Breaking]{Solutions of the Camassa--Holm equation with accumulating 
breaking times}

\author[K. Grunert]{Katrin Grunert}
\address{Department of Mathematical Sciences\\ NTNU\\Norwegian University of Science 
and Technology\\ NO-7491 Trondheim\\ Norway}
\email{\mailto{katring@math.ntnu.no}}
\urladdr{\url{http://www.math.ntnu.no/~katring/}}

\thanks{Research supported by the grant {\it Waves and Nonlinear Phenomena (WaNP)} from the Research Council of Norway.}

\begin{abstract}
We present two initial profiles to the Camassa--Holm equation which yield 
solutions with accumulating breaking times.  
\end{abstract}

\maketitle

\section{Introduction}

The Camassa--Holm (CH) equation \cite{CH, CH2}
\begin{equation}
u_t-u_{txx}+3uu_x-2u_xu_{xx}-uu_{xxx}=0,
\end{equation}
which serves as a model for shallow water \cite{CL},
has been studied intensively over the last twenty years, due to its rich 
mathematical structure. For example, it is bi-Hamiltonian \cite{FF}, completely integrable \cite{C2}
and has infinitely many conserved quantities, see e.g. \cite{L3}. Yet another property attracted 
considerable attention: Even smooth initial data can lead to classical 
solutions, which only exist locally due to wave breaking, see e.g. \cite{CE, 
CE2, CE3}. That means the spatial derivative $u_x(t,.)$ of the solution $u(t,.)$ 
becomes unbounded within finite time, while $\norm{u(t,.)}_{H^1}$ remains bounded. In 
addition, energy concentrates on sets of measure zero when wave breaking takes 
place. Neglecting this concentration yields to a dissipation and hence to the 
so-called dissipative solutions \cite{BC2, HR2}. However, taking care of the 
energy, yields another class of solutions, the conservative ones \cite{BC,HR}. 
Moreover, it is also possible to take out only a fraction of the concentrated 
energy, giving rise to the recently introduced $\alpha$-dissipative solutions 
\cite{GHR6}. A very illustrating example for this phenomenon is given by the 
so-called peakon-antipeakon solutions, which enjoy wave breaking and therefore 
can be prolonged thereafter in various ways as presented in detail in 
\cite{GH,GHR6}.

However, as already the study of the peakon-antipeakon solutions shows, there 
are only very few solutions of the CH equation, which can be computed 
explicitly. Even in the case of travelling wave solutions, which have been 
classified by J. Lenells in \cite{L1}, some of them are only given implicitely. 
Having a close look at the construction of various types of solutions \cite{BC, 
BC2, GHR, HR, HR2}, reveals that they are based on a reformulation of the CH 
equation as a system of ordinary differential equations in a suitable Banach 
space by a generalized method of characteristics. Thus computing solutions 
explicitly, would involve a change of variables from Eulerian to Lagrangian 
coordinates, solving the equation in Lagrangian coordinates and finally changing 
back from Lagrangian to Eulerian coordinates, as outlined in 
Section~\ref{background}. A task which in many cases seems to be impossible. 
Thus most results concerning the prediction of wave breaking are obtained by 
following solutions along characteristics, see \cite{C,CE,G}.

However, a good understanding of solutions along characteristics allows the 
prediction of wave breaking in the nearby future, and in the case of 
conservative solutions, which one can follow both forward and backward in time, 
also to find out if wave breaking occurred recently, \cite{G}. 
Based on this knowledge it is possible to construct some initial data $u_0(x)$, which has $t=0$ as an accumulation point of breaking times in the conservative case, as we will see in Section~\ref{sec:ex}. 
We will prove the following result.

\begin{theorem}\label{thm:accu}
Let $q\in(0,1)$ and 
\begin{equation}
u_0(x)=\begin{cases} 
 \frac12 q(1+\frac{1-q^4}{2}x), &\quad \text{ for }x\in[-\frac{2}{1-q^4},0],\\
k_jx+d_j, & \quad \text{ for } x\in[x_j,x_{j+1}], j\in\mathbb{N}\\
 0, & \quad \text{ otherwise},
\end{cases}
\end{equation}
where the endpoints of the intervals $[x_j, x_{j+1}]$ are inductively defined through
\begin{equation}
 x_0=0, \quad  x_{2j+2}-x_{2j+1}=x_{2j+1}-x_{2j}=q^{4j}
\end{equation}
and the constants $k_j$ and $d_j$ satisfy
\begin{align*}
k_{2j}&  =-\frac{1}{q^{j-1}},   \quad k_{2j+1}=\frac12\frac{q+q^4}{q^j}, \\
d_{2j}& =\frac12 \frac{1}{q^{j-1}(1-q^4)} (4-3q^{4j}-q^{4(j+1)}), \\
d_{2j+1}&=-\frac12 \frac{1}{q^{j-1}(1-q^4)}(2+2q^3-q^{4j+3}-2q^{4j+4}-q^{4j+7}).
\end{align*}
Furthermore, denote by $u(t,x)$ the conservative solution of the CH equation with initial data $u(0,x)=u_0(x)$, then $t=0$ is an accumulation point of breaking times.
\end{theorem}

In Section~\ref{sec:cusp} we are going to have a closer look at the cuspons with exponential decay \cite{L1}, a class of travelling wave solutions with non-vanishing asymptotics, for which wave breaking takes place at any time. 
However, for any fixed time $t$, the Radon measure is purely absolutely continuous, which means that no energy is concentrated on sets of measure zero. This means in particular, that the set of points where wave breaking occurs at time $t$, consists of a single point both in Eulerian and Lagrangian coordinates, as we will see. Additional we are going to show that the breaking point is not traveling along one characteristic with respect to time, but is metaphorically speaking, jumping from one characteristic to the next one. 
These observations are very interesting since usually wave breaking is linked to the concentration of energy on sets of measure zero in Eulerian coordinates, which corresponds to wave breaking taking place on sets of positive measure in Lagrangian coordinates. Thus it is natural that manipulating the concentrated energy gives rise to different solution concepts. For this example however the question turns up of how the dissipative solution looks like? Does it coincide with the conservative one or not?

Finally, Section~\ref{Observation} concludes this note, by an 
observation concerning the concentration of energy in the case of accumulating 
breaking times.

\section{Background material}\label{background}

The aim of this section is to outline the construction of conservative solutions 
of the CH equation, which consists of two main parts. On the one hand the 
interplay of Eulerian and Lagrangian coordinates and on the other hand the 
reformulation of the CH equation in Lagrangian coordinates. We will restrict 
ourselves to presenting those results, which are going to play a key role in 
what follows. For details we refer the interested reader to \cite{G} and 
\cite{GHR6}.

Let us start with the interplay between Eulerian and Lagrangian coordinates. It 
is well-known that solutions of the CH equation might enjoy wave breaking within 
finite time. This means that the solution itself remains bounded while its 
spatial derivative becomes unbounded from below pointwise. In particular, energy 
concentrates on sets of measure zero at breaking time, which can be described 
with the help of positive, finite Radon measures. Thus, the admissible set of 
Eulerian coordinates, $\D$, is defined as follows.

\begin{definition} [Eulerian coordinates]\label{def:D} The set $\D$ is
  composed of all pairs $(u,\mu)$ such that
  $u\in H^1(\Real)$ and $\mu$ is a positive, finite
  Radon measure whose absolutely continuous part,
  $\muac$, satisfies
\begin{equation}
\label{eq:abspart}
\muac=u_x^2\,dx.
\end{equation}
\end{definition}

Rewriting the CH equation in the weak formulation yields
\begin{equation}\label{weakCH}
u_t+uu_x+P_x=0
\end{equation}
where 
\begin{equation}
P(t,x)=\frac12 \int_\Real e^{-\vert x-z\vert}u^2(t,z)dz+\frac14 \int_\Real 
e^{-\vert x-z\vert}d\mu(t,z).
\end{equation}
A close inspection of \eqref{weakCH} reveals that one can try to compute 
solutions of the CH equation using the method of characteristics. Indeed, this 
is possible but only under the assumption that $\mu$ is a purely absolutely 
continuous Radon measure and that we are given some initial characteristic 
$y_0(\xi)=y(0,\xi)$. Moreover, due to wave breaking we can only expect to obtain 
local solutions so far. 
Thus $y(t,\xi)$ is the solution to  
\begin{equation}\label{eq:char}
 y_t(t,\xi)=u(t,y(t,\xi)),
\end{equation}
for some given initial data $y_0(\xi)=y(0,\xi)$,
and additional we introduce the function 
\begin{equation}\label{eq:map}
 U(t,\xi)=u(t,y(t,\xi)) , 
\end{equation}
whose time evolution is given through \eqref{weakCH}. 

In the general case where $\mu$ is a positive finite Radon measure, we use the 
following mapping from Eulerian coordinates $\D$ to Lagrangian coordinates $\F$ 
to obtain an admissible initial characteristic $y_0(\xi)$ for any initial data 
$(u_0,\mu_0)\in\D$.

\begin{definition} 
\label{th:Ldef}
For any $(u,\mu)$ in $\D$, let
\begin{subequations}
\label{eq:Ldef}
\begin{align}
\label{eq:Ldef1}
y(\xi)&=\sup\left\{y\ |\ \mu((-\infty,y))+y<\xi\right\},\\
\label{eq:Ldef2}
h(\xi)&=1-y_\xi(\xi),\\
\label{eq:Ldef3}
U(\xi)&=u\circ{y(\xi)}.
\end{align}
\end{subequations}
Then $(y,U,h)\in\F$. We denote by $L\colon \D\rightarrow \F$ the mapping which 
to any element $(u,\mu)\in\D$ associates $X=(y,U,h)\in \F$ given by 
\eqref{eq:Ldef}.  
\end{definition}

The big advantage of this change of variables is due to the fact that $y$, $U$, 
and $h$ are all functions, and they will remain functions for all times. However 
there are a lot of properties these functions have to satisfy to enable us to 
construct global (conservative) solutions to the CH equation. All of them are 
collected in the following definition.

\begin{definition}[Lagrangian coordinates]\label{def:F} The set $\F$ is composed 
of all $X=(\zeta,U,h)$, such that 
 \begin{subequations}
\label{eq:lagcoord}
\begin{align}
\label{eq:lagcoord1}
& (\zeta, U,h,\zeta_\xi,U_\xi)\in L^\infty(\Real)\times [L^\infty(\Real)\cap 
L^2(\Real)]^4,\\
\label{eq:lagcoord2}
&y_\xi\geq0, h\geq0, y_\xi+h>0
\text{  almost everywhere},\\
\label{eq:lagcoord3}
&y_\xi h=U_\xi^2\text{ almost everywhere},\\
\label{eq:lagcoord4}
&y+H\in G,
\end{align}
\end{subequations}
where we denote $y(\xi)=\zeta(\xi)+\xi$ and $H(\xi)=\int_{-\infty}^\xi 
h(\eta)d\eta$. 
\end{definition}

 Here we denote by $G$ the subgroup of the group of homeomorphisms from $\Real$ 
to $\Real$ such that 
\begin{subequations}
\label{eq:Gcond}
 \begin{align}
  \label{eq:Gcond1}
  f-\id \text{ and } f^{-1}-\id &\text{ both belong to } W^{1,\infty}(\Real), \\
  \label{eq:Gcond2}
  f_\xi-1 &\text{ belongs to } L^2(\Real),
 \end{align}
 \end{subequations}
where $\id$ denotes the identity function.
In particular, $G$ coincides with the set of relabelling functions, which enable 
us to identify equivalence classes in $\F$. This is necessary since we have 3 
Lagrangian coordinates in contrast to 2 Eulerian coordinates. 

In the case of conservative solutions, the reformulation of the CH equation in 
Lagrangian coordinates is given through
\begin{subequations}
  \label{eq:chsyseq}
  \begin{align}
    \label{eq:chsyseq1}
    \zeta_t&=U,\\
    \label{eq:chsyseq:2}
    \zeta_{\xi,t}& =U_\xi,\\
    \label{eq:chsyseq3}
    U_t&=-Q,\\
    \label{eq:chsyseq4}
    U_{\xi,t}& =\frac{1}{2}h+(U^2-P)y_\xi, \\
    \label{eq:chsyseq5}
    h_t&=2(U^2-P)U_\xi,
  \end{align}
\end{subequations}
where 
\begin{equation}
 \label{eq:P}
P(t,\xi)=\frac14\int_\Real
  e^{-\abs{y(t,\xi)-y(t,\eta)}}(2U^2y_\xi+h)(t,\eta)\,d\eta,
\end{equation}
and 
\begin{equation}
 \label{eq:Q}
Q(t,\xi)=-\frac14\int_\Real \sign(\xi-\eta)e^{-\abs{y(t,\xi)-y(t,\eta)}}(2
  U^2y_\xi+h)(t,\eta)\,d\eta.
\end{equation}
One can show that both $P(t,.)$ and $Q(t,.)$ belong to $H^1(\Real)$ and that this system of ordinary differential equations admits global unique solutions in 
$\F$. 

Hence it remains to get back from Lagrangian to Eulerian coordinates 
\begin{definition}
\label{th:umudef} 
Given any element $X=(y,U,h)\in\F$, we define $(u,\mu)$ as follows
\begin{subequations}
\label{eq:umudef}
\begin{align}
\label{eq:umudef1}
u(x)&=U(\xi) \text{ for any }\xi \text{ such that }x=y(\xi),\\
\label{eq:umudef2}
\mu&=y_\#(h(\xi)\,d\xi),
\end{align}
\end{subequations}
We have that $(u,\mu)$ belongs to $\D$. We
denote by $M\colon \F\rightarrow\D$ the mapping
which to any $X$ in $\F$ associates the element $(u,\mu)\in \D$ as given
by \eqref{eq:umudef}. 
\end{definition}
In particular, the mapping $M$ maps elements $X\in\F$ belonging to one and the 
same equivalence class in $\F$ to one and the same element in $\D$, see 
\cite{HR}.

Before focusing on results related to our further investigations, we would like 
to emphasize that there are several ways to prolong solutions past wave breaking. 
Most of them are related to how the energy is manipulated at breaking time. In 
the case of conservative solutions we mean that the total amount of energy 
remains unchanged with respect to time, that is 
\begin{equation*}
\norm{u(t,.)}_{L^2}^2+\mu(t,\Real)=\norm{u_0}_{L^2}^2+\mu_0(\Real).
\end{equation*} 
For a detailed discussion of this topic we refer to \cite{GHR6} and for more 
information on how wave breaking is described in Eulerian and Lagrangian 
coordinates we refer to \cite{G}.

Of great importance for us, will be the prediction of wave breaking. Although we 
cannot determine exactly at which points $(t,x)$ wave breaking takes place, we 
can at least determine if wave breaking occurs in the nearby future or not. The 
result is contained in the following theorem, which is a slight reformulation of 
\cite[Theorem 1.1]{G}.
\begin{theorem}\label{thm:main}
 Given $(u_0,\mu_0)\in \D$ and denote by $(u(t),\mu(t))\in\D$ the global 
conservative solution of the CH equation at time $t$ with initial data $(u_0,\mu_0)$. Moreover, let 
$C=2(\norm{u_0}_{L^2}^2+\mu_0(\Real))$, then the following holds:
\begin{enumerate}
 \item If $u_{0,x}(x)< -\sqrt{2C}$ for some $x\in\Real$, then wave breaking will 
occur within the time interval $[0,T]$, where $T$ denotes the solution of 
 \begin{equation}
 \frac{u_{0,x}(x)+\sqrt{2C}}{u_{0,x}(x)-\sqrt{2C}}= \exp(-\sqrt{2C}T).
\end{equation}
\item If $u_{0,x}(x)> \sqrt{2C}$ for some $x\in\Real$, then wave breaking 
occured within the time interval $[T,0]$, where $T$ denotes the solution of 
 \begin{equation}
 \frac{u_{0,x}(x)+\sqrt{2C}}{u_{0,x}(x)-\sqrt{2C}}= \exp(-\sqrt{2C}T).
\end{equation}
\end{enumerate}
\end{theorem}

As far as the situation in Lagrangian coordinates is concerned much more is 
known, due to the fact that the prediction of wave breaking is based on 
following solutions along characteristics.
A key role, in that context, plays the following set,
\begin{equation}\label{eq:defKgamma}
  \kappa_{1-\gamma}=\{\xi\in\Real\mid \frac{h_0}{y_{0,\xi}+h_0}(\xi)\geq 
1-\gamma\text{, } U_{0,\xi}(\xi)\leq 0\}, \quad \gamma\in [0,\frac12].
\end{equation}
Every condition imposed on points $\xi\in\kappa_{1-\gamma}$ is motivated by what 
is known about wave breaking. Indeed, if wave breaking occurs at some time 
$t_b$, then energy is concentrated on sets of measure zero in Eulerian 
coordinates, which correspond to the sets where $\frac{h}{y_\xi+h}(t_b,\xi)=1$ 
in Lagrangian coordinates. Furthermore, at time $t_b$ the solution $u$ is 
bounded while its spatial derivative $u_x$ becomes unbounded from below 
pointwise, see \cite{ConstantinIvanov:2008, GY}. In Lagrangian coordinates this 
means that $U_\xi(t_b,\xi)=0$ and $U_\xi(t,\xi)$ changes sign from negative to 
positive at breaking time, \cite{G}. 

The next lemma, which is a reformulation of \cite[Corollary 18]{GHR6}, enables 
us to predict wave breaking in the nearby future along characteristics.

\begin{lemma}\label{lem:1}
Let $X_0\in\F$. Denote by $X=(\zeta, U, \zeta_\xi, U_\xi, h)\in C(\Real_+,\F)$ 
the global solution of \eqref{eq:chsyseq} with initial data $X_0$ and by 
$\tau_1(\xi)\geq 0$ the first time wave breaking occurs at the point $\xi$. 
Moreover, let $M=\norm{U_0^2y_{0,\xi}+h_0}_{L^1}$. Then the following statements 
hold:
  
  (\textit{i}) We have 
  \begin{equation}\label{eq:G3sol}
    \norm{\frac{1}{y_\xi+h}(t,\dott)}_{L^\infty}\leq 
2e^{C(M)T}\norm{\frac{1}{y_{0,\xi}+h_0}}_{L^\infty},
  \end{equation}
  and 
  \begin{equation}
    \label{eq:G3bsol}
    \norm{(y_\xi+h)(t,\dott)}_{L^\infty}\leq 
2e^{C(M)T}\norm{y_{0,\xi}+h_0}_{L^\infty}
  \end{equation}
  for all $t\in[0,T]$ and a constant $C(M)$ which depends on $M$. 

  (\textit{ii}) There exists a $\gamma\in(0,\frac12)$ depending only on $M$ such 
that if $\xi\in \kappa_{1-\gamma}$, then $\frac{y_\xi}{y_\xi+h}(t,\xi)$ is a 
decreasing function and $\frac{U_\xi}{y_\xi+h}(t,\xi)$ is an increasing 
function, both with respect to time for $t\in[0,\min(\tau_1(\xi),T)]$.  
Therefore we have 
  \begin{equation}\label{eq:G6sol}
    \frac{U_{0,\xi}}{y_{0,\xi}+h_0}(\xi)\leq\frac{U_\xi}{y_\xi+h}(t,\xi)\leq 
0\quad \text{and} \quad 0\leq \frac{y_\xi}{y_\xi+h}(t,\xi)\leq 
\frac{y_{0,\xi}}{y_{0,\xi}+h_0}(\xi),
  \end{equation}
  for $t\in[0,\min(\tau_1(\xi),T)]$.
  In addition, for $\gamma$ sufficiently small, depending only on $M$ and $T$, 
we have 
  \begin{equation}\label{eq:G5sol}
    \kappa_{1-\gamma}\subset \{\xi\in\Real\mid 0\leq \tau_1(\xi)<T\}.
  \end{equation}

  (\textit{iii})
  For any given $\gamma>0$, there exists $\hat T>0$ such that 
  \begin{equation}\label{eq:G4sol}
    \{\xi\in\Real\mid 0<\tau_1(\xi)<\hat T\}\subset \kappa_{1-\gamma}.
  \end{equation}
\end{lemma}

$(\textit{i})$ ensures that the function $\frac{h}{y_\xi+h}(t,\xi)$ is 
well-defined. $(\textit{ii})$ gives us a possibility to fix at first some time 
interval $[0,T]$ and thereafter by finding a suitable $\gamma$ to identify 
points which enjoy wave breaking for sure. $(\textit{iii})$ on the other hand 
gives us a possibility to choose $\gamma$ first and identifying a time interval 
$[0,T]$ thereafter, such that $\kappa_{1-\gamma}$ contains all points enjoying 
wave breaking within $[0,T]$. 

Obviously the question occurs if a point $\xi$ in Lagrangian coordinates can 
enjoy wave breaking infinitely many times within a fixed time interval $[0,T]$. 
According to \cite[Corollary 19]{GHR6}, which we state here for the sake of 
completeness, this is not possible.

\begin{lemma}[\cite{GHR6} Corollary 19]\label{lem:2}
  Denote by $X(t)=(y,U,y_\xi, U_\xi,h)(t)$ the global solution of 
\eqref{eq:chsyseq} with $X(0)=X_0\in\F$ in $C(\Real_+,\F)$. Let $M=\norm{ 
U_0^2y_{0,\xi}+h_0}_{L^1}$ and denote by $\tau_j(\xi)\geq 0$ the j'th time wave 
breaking occurs at a point $\xi\in\Real$. Then for any $\xi\in\Real$ the 
possibly infinite sequence
  $\tau_j(\xi)$ cannot accumulate.

In particular, there exists a time $\hat T$ depending on $M$ such that any point 
$\xi$ can experience
  wave breaking at most once within the time interval $[T_0,T_0+\hat T]$ for any 
$T_0\ge0$. More
  precisely, given $\xi\in\Real$, we have
  \begin{equation}
  \tau_{j+1}(\xi)-\tau_j(\xi)>\hat T \text{ for all } j. 
  \end{equation}
  In addition,
  for $\hat T$ sufficiently small, we get that in this case $U_\xi(t,\xi)\geq 0$ 
for all
  $t\in[\tau_j(\xi), \tau_j(\xi)+\hat T]$.
\end{lemma}

Here it is important to note that we can only say that the sequence 
$\tau_j(\xi)$ for $\xi\in\Real$ does not accumulate in Lagrangian coordinates. 
We have, however, no possibility to conclude that the same result holds in 
Eulerian coordinates for $x\in\Real$, since we follow the solution in Lagrangian 
coordinates along characteristics. 

\begin{remark}
 In Section~\ref{sec:cusp} we are going to look at the case of a traveling wave solution $u(t,x)$ with nonvanishing, but equal asymptotics as $x\to\pm\infty$. That is, there exists $c\in\Real$ such that $u(t,x)-c\in H^1(\Real)$ for all $t\in\Real$. 
 Also in this case the description of solutions in Lagrangian coordinates is possible, by slightly changing the definition of $\D$ and $\F$, while leaving Definition~\ref{th:Ldef}, Definition~\ref{th:umudef} and \eqref{eq:chsyseq}--\eqref{eq:Q} unchanged. To be more explicit we have to replace $U(\xi)\in H^1(\Real)$ by $U(\xi)-c\in H^1(\Real)$ in Definition~\ref{def:F}, while the remaining assumptions remain unchanged. Moreover, $P(t,\xi)-c^2$ and $Q(t,\xi)$ belong to $H^1(\Real)$ for all $t\in\Real$. 
 
The set of points in Lagrangian coordinates where wave breaking occurs at time $t$ still coincides with 
\begin{equation*}
 \{ \xi\in\Real \mid y_\xi(t,\xi)=0\}.
\end{equation*} For details we refer the interested reader to \cite{GHR2} and \cite{GHR}. 
\end{remark}

\section{Can breaking times accumulate for solutions of the CH equation?}\label{sec:ex}
The aim of this section is to identify some solutions of the CH equation with an accumulating sequence of breaking times. Since it is nearly impossible, beside of some special cases, to compute solutions explicitly, we constructed some initial data, whose corresponding conservative solution has an accumulating sequence of breaking times. The construction is based on Lemma~\ref{lem:1} and Lemma~\ref{lem:2}, which make it unnecessary to compute the actual solution.

To be more specific, we aimed at finding some initial data $(u_0, 
\mu_0)\in \D$, with $\mu_0=\mu_{0,ac}$ such that 
\begin{itemize}
\item $u_0$ has compact support, that is $supp(u_0)=[-c, c]$ for some $c>0$,
\item $u_0$ is a piecewise linear and continuous function, that is there exists 
an increasing sequence $x_j$ and two sequences $k_j$ and 
$d_j$ such that
\begin{enumerate}
\item $x_{-1}=-c$ and $x_0=0$
\item $(x_{j+1}-x_j)\to 0$ as $j\to \infty$ (nonincreasing)
\item $k_{2j}<0$ and $k_{2j+1}>0$ for $j\in\mathbb{N}$
\item $ -k_{2j} \to \infty$ and $k_{2j+1}\to\infty$ as $j\to\infty$ (strictly 
increasing)
\item $u_0(x)=k_jx+d_j$ for $x\in[x_j, x_{j+1}]$, $j=-1,0,1,2,\dots$
\end{enumerate}
\item $u_0\in H^1(\Real)$
\end{itemize}

Each of these assumptions is motivated by what is known about the prediction of wave breaking. Hence we want to have a close look at them, before turning to the proof of Theorem~\ref{thm:accu}.
Although $u_0$ has compact support, $u_{0,x}$ is not going to be bounded on 
$[0,c]$. Thus in order for $u_0$ to be in $H^1(\Real)$ there must be a balance 
between the increasing sequences $-k_{2j}$ and $k_{2j+1}$ and the decreasing 
sequence $x_{j+1}-x_j$. Moreover, the condition $\vert k_j\vert\to\infty$ is 
necessary to impose since Lemma~\ref{lem:1} $(\textit{ii})$ points out that 
given some time interval $[0,T]$, we can be sure that wave breaking occurs for 
all points which lie inside $\kappa_{1-\gamma}$ for $\gamma$ small enough. Thus 
if we can find for each $\gamma>0$ infinitely many nonintersecting intervals in 
Lagrangian coordinates which lie inside $\kappa_{1-\gamma}$, the claim follows. 
So what does $\xi\in \kappa_{1-\gamma}$ mean? If $\mu_0$ is absolutely 
continuous, which we assume, then 
\begin{equation}
y_{0,\xi}(\xi)=\frac{1}{1+u_{0,x}^2(y_0(\xi))}=\frac{1}{1+k_j^2} , \quad 
\text{for } y_0(\xi)\in[x_j,x_{j+1}],
\end{equation}
which means that $y_0(\xi)$ is an increasing function. 
Thus 
\begin{equation}\label{eq:1}
\frac{h_0(\xi)}{y_{0,\xi}(\xi)+h_0(\xi)}=1-\frac{y_{0,\xi}(\xi)}{y_{0,\xi}
(\xi)+h_0(\xi)}=1-y_{0,\xi}(\xi)=1-\frac{1}{1+k_j^2},
\end{equation}
for $y_0(\xi)\in[x_j,x_{j+1}]$, where we used $y_{0,\xi}(\xi)+h_0(\xi)=1$.
Furthermore,
\begin{equation}\label{eq:2}
U_{0,\xi}(\xi)=u_{0,x}(y_0(\xi))y_{0,\xi}(\xi)=\frac{k_j}{1+k_j^2}, \quad 
\text{for } y_0(\xi)\in[x_j,x_{j+1}].
\end{equation}
Combining \eqref{eq:1} and \eqref{eq:2} yields $\xi\in \kappa_{1-\gamma}$ if and 
only if $y_0(\xi) \in [x_j,x_{j+1}]$ with $k_j<0$ and $\gamma\geq 
\frac{1}{1+k_j^2}$. Thus our assumptions on $k_j$ are chosen in such a way that 
$\kappa_{1-\gamma}$ consists of infinitely many nonintersecting intervals.

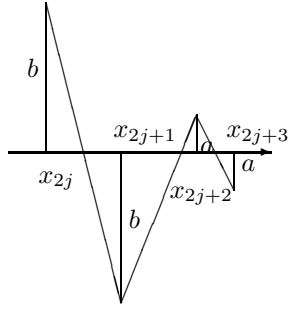
\begin{figure}\centering
\begin{picture}(3.5,5)
\put(0,2.5){\vector(1, 0){3.5} }
\put(0.5,4.5){\line(0,-1) {2} }
\put(1.5,2,5){\line(0,-1) {2} }
\put(0.5,4.5){\line(1,-4) {1} }
\put(2.5,2.5){\line(0,1){0.5}}
\put(1.5,0.5){\line(2,5) {1}}
\put(3,2.5){\line(0,-1){0.5}}
\put(2.5,3){\line(1,-2){0.5}}
\put(0.25,3.5){$b$}
\put(1.6,1.5){$b$}
\put(3.1,2.25){$a$}
\put(2.55,2.5){$a$}
\put(0.4,2.1){$x_{2j}$}
\put(1.4, 2.7){$x_{2j+1}$}
\put(2.15,1.9){$x_{2j+2}$}
\put(2.9, 2.7){$x_{2j+3}$}
\end{picture}
\caption{An illustration of the piecewise linear function 
$u$}\label{fig:construe}
\end{figure}

\begin{proof}[Proof of Theorem 1.1]
That $u_0\in H^1(\Real)$ relies on the fact that the geometric series $\sum_{n=0}^\infty q^n$ converges for $q\in (0,1)$ and the proof is left to the reader. We rather focus on showing the existence of a convergent sequence of breaking times.

Let $M=\norm{u_0}_{H^1}^2=\norm{U_0^2y_{0,\xi}+h_0}_{L^1}$, then according to 
Lemma~\ref{lem:2} there exists a time $\hat T$ depending on $M$ such that any 
point $\xi\in\Real$ in Lagrangian coordinates can experience wave breaking at 
most once within the time interval $[0,\hat T]$. Moreover, Lemma~\ref{lem:1} 
(ii) implies that there exists $\hat\gamma\in (0,\frac12)$ depending only on $M$ 
and $\hat T$ such that 
\begin{equation}
\kappa_{1-\hat\gamma}\subset\{\xi\in\Real\mid 0\leq \tau_1(\xi)< \hat T\},
\end{equation}
and such that \eqref{eq:G6sol} holds for all $\xi\in\kappa_{1-\hat\gamma}$.  

Let $\gamma\in(0,\hat\gamma)$ and $\xi\in\kappa_{1-\gamma}$, then we have 
\begin{equation}\label{eq:derivative}
\left(\frac{U_\xi}{y_\xi+h}\right)_t=\frac12 +\left( 
U^2-P-\frac12\right)\frac{y_\xi}{y_\xi+h}-(2U^2-2P+1)\frac{U_\xi^2}{(y_\xi+h)^2}
,
\end{equation}
and in particular $\frac{U_\xi}{y_\xi+h}(t,\xi)$ changes sign from negative to 
positive at breaking time $\tau_1(\xi)<\hat T$. Since 
$\norm{(U^2-P)(t,.)}_{L^\infty}$ can be bounded uniformly by a constant 
depending on $M$, we have that 
\begin{equation}\label{eq:updown}
\frac12-C(M)\gamma\leq 
\left(\frac{U_\xi}{y_\xi+h}\right)_t\leq\frac12+C(M)\gamma,
\end{equation}
for some constant $C(M)$ only dependent on $M$. Assume additional that $\gamma$ 
is so small that $\frac12-C(M)\gamma>0$, then \eqref{eq:updown} allows us to 
derive an upper and a lower bound on $\tau_1(\xi)$ for $\xi\in 
\kappa_{1-\gamma}$ such that $\frac{h_0}{y_{0,\xi}+h_0}(\xi)=1-\gamma$. Indeed, 
on the one hand \eqref{eq:updown} implies that 
\begin{equation}
\left(\frac12 -C(M)\gamma\right) t-\sqrt{\gamma(1-\gamma)}\leq 
\frac{U_\xi}{y_\xi+h}(t,\xi),
\end{equation}
and hence $\tau_1(\xi)\leq T_{1,\gamma}=\frac{\sqrt{\gamma(1-\gamma)}}{\frac12 
-C(M)\gamma}$.
On the other hand, we have 
\begin{equation}
\frac{U_\xi}{y_\xi+h}(t,\xi)\leq \left(\frac12+C(M)\gamma 
\right)t-\sqrt{\gamma(1-\gamma)},
\end{equation}
which implies 
$T_{2,\gamma}=\frac{\sqrt{\gamma(1-\gamma)}}{\frac12+C(M)\gamma}\leq 
\tau_1(\xi)$.
Thus $T_{2,\gamma}\leq \tau_1(\xi)\leq T_{1,\gamma}$ for all 
$\xi\in\kappa_{1-\gamma}$ such that $\frac{h_0}{y_{0,\xi}+h_0}(\xi)=1-\gamma$. 
Seen as a function of $\gamma$, both $T_{1,\gamma}$ and $T_{2,\gamma}$ are 
strictly decreasing or equivalently 
\begin{equation}
\lim_{\gamma\to 0} T_{1,\gamma}=0 \quad \text{ and }\quad \lim_{\gamma\to 0} 
T_{2,\gamma}=0.
\end{equation}
By definition the sequence $k_{2j}^2=\frac{1}{q^{2(j-1)}}\to \infty$ as $j\to\infty$. Thus the corresponding
sequence $\gamma_j=\frac{1}{1+\kappa_{2j}^2}\to 0$ as $j\to\infty$ and each $\kappa_{1-\gamma_j}$ consists of infinitely many non-intersecting intervals. Thus choosing to every $j\in\mathbb{N}$ a point $z_j\in[x_{2j}, 
x_{2j+1}]$ with breaking time $t_j$, the above argument shows that $t_j\to 0$.
\end{proof}

\begin{remark}
Before continuing with another interesting example, we would like to point out 
another interesting fact about our constructed initial data $u_0$. According to 
Theorem~\ref{thm:main} it is not only possible to predict if wave breaking 
occurs in the nearby future or not, but also if wave breaking occurred recently 
or not. Moreover, it is also possible to adapt in the conservative case 
Lemma~\ref{lem:1} and Lemma~\ref{lem:2} to going backward in time. Thus 
following a similar argument to the one presented above one can also show that 
there exists an increasing sequence of breaking times with limiting value zero. 
\end{remark}

\section{Cuspons with exponential decay - from a wave breaking point of view}\label{sec:cusp}

In \cite{L1, L2} J. Lenells classified all travelling wave solutions of the 
Camassa--Holm equation. One of them, namely the cuspons with exponential decay, 
serve as an interesting example of solutions enjoying wave breaking. On the one 
hand wave breaking occurs at any time and on the other hand the energy is not 
concentrated on sets with positive measure in Lagrangian coordinates, but in a 
single point which does not correspond to the singular continuous part of the 
positive, finite Radon measure.

\begin{definition}[Cuspon with exponential decay]\label{def:cusp}
Given $m$, $s$, and $M$ such that $m<s<M$, let $\kappa=\frac12(s-2m-M)$, then 
the cuspon with exponential decay and speed $s+\kappa$ is defined through 
\[
u(t,x)=\phi(x-(s+\kappa)t)+\kappa
\]
where $\phi(x)$ is implicitly given through
\begin{equation}\label{def:phi}
\phi_x^2=\frac{(M-\phi)(\phi-m)^2}{(s-\phi)}
\end{equation}
and satisfies 
\begin{subequations}\label{cond:phi}
\begin{align}
& \phi(-x)=\phi(x) \\
& \phi(0)=s\\
& \phi_x(x)<0 \quad \text{for} \quad x>0\\
& \lim_{x\to\infty} \phi(x)=m.
\end{align}
\end{subequations}
\end{definition}

\begin{figure}\centering
\includegraphics[width=5cm]{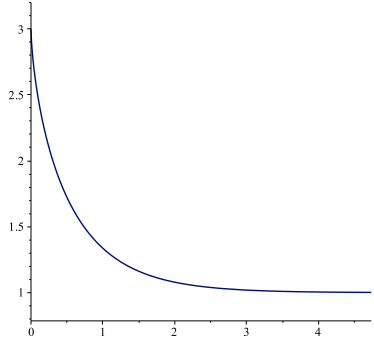}
\caption{Plot of the function $\phi(x)$ for $m=1$, $s=3$, and $M=5$.}
\end{figure}

Note that for any $\phi(x)$ defined through \eqref{def:phi} and 
\eqref{cond:phi}, we have that $m\leq \phi(x)\leq s$ for all $x\in\Real$. Thus 
$\phi_x(x)$ is finite for all $x\in\Real\backslash\{0\}$, 
\begin{equation*}
\lim_{x\to 0-} \phi_x(x)=\infty \quad \text{ and } \quad \lim_{x\to 0+} 
\phi_x(x)=-\infty
\end{equation*}
and in particular 
\begin{equation*}
\lim_{x\to 0} \phi_x^2(x)=\infty.
\end{equation*}
This implies for $u(t,x)=\phi(x-(s+\kappa)t)+\kappa$, that $u_x(t,x)$ is 
well-defined for all $x\in\Real\backslash\{(s+\kappa)t\}$ and 
\[
\lim_{x\to (s+\kappa)t} u_x^2(t,x)=\infty \quad \text{ for all } t.
\]
Or in order words wave breaking occurs for the solution $u(t,x)$ defined in 
Definition~\ref{def:cusp} for all points $(t,x)$ such that $x=(s+\kappa)t$.

{\bf Breaking points in Lagrangian coordinates.}
Let $u(0,x)=\phi(x)+\kappa$ and denote by $y_0(\xi)=y(0,\xi)$ the initial 
characteristics in Lagrangian coordinates given through \eqref{eq:Ldef1}. Then 
we claim that the point $x=0$ in Eulerian coordinates is mapped to the single 
point $\bar \xi$ in Lagrangian coordinates which satisfies  $y_{0}(\bar \xi)=0$ 
and  
\begin{equation}\label{break:set}
\{ \bar \xi \}=\{ \xi\in\Real\mid y_{0,\xi}(\xi)=0\}.
\end{equation}

Let us start by considering the function $g(x)=x+\int_{-\infty}^x 
\phi^2_x(y)dy$, which has as a pseudo inverse the function $y(\xi)$ (cf. 
\eqref{eq:Ldef1}). If $\int_{-\infty}^x\phi^2_x(y)dy$ exists 
for all $x\in(-\infty,0]$, then it follows from the symmetry of $\phi(x)$ that 
the function $g(x)$ is well-defined, strictly increasing and bijektiv. This 
means in particular to any $x\in\Real$ there exists a unique $\xi\in\Real$ such 
that $y(\xi)=x$. Indeed, we have for $x\leq 0$,
\begin{align*}
\int_{-\infty}^{x} \phi_x^2(y)dy & =\int_{-\infty}^x 
\frac{\sqrt{M-\phi(y)}(\phi(y)-m)}{\sqrt{s-\phi(y)}}\phi_x(y)dy\\
& = \int_m^{\phi(x)} \frac{\sqrt{M-z}(z-m)}{\sqrt{s-z}}dz\\
& \leq \sqrt{M-m}(s-m)\int_m^s \frac{1}{\sqrt{s-z}}dz\\
& = 2\sqrt{(M-m)(s-m)}(s-m)<\infty
\end{align*}
and hence $y(\xi)$ is bijective.

Denote by $\bar \xi$ the unique $\xi\in\Real$ such that $y(\bar \xi)=0$. Then we 
have that $\phi_x(x)$ is well-defined for all $x\in\Real\backslash \{0\}$ and 
hence 
\[
y_{0,\xi}(\xi)=\frac{1}{1+\phi_x^2(y_0(\xi))}=\frac{(s-\phi(y_0(\xi)))}{
(s-\phi(y_0(\xi)))+(M-\phi(y_0(\xi)))(\phi(y_0(\xi))-m)^2}
\] for all $\xi\in\Real\backslash\{\bar \xi\}$. In particular, 
$y_{0,\xi}(\xi)>0$ for all $\xi\in \Real\backslash \{\bar \xi\}$. Thus if we can 
show that $y_{0,\xi}(\bar \xi)=0$, we obtain as an immediate consequence \eqref{break:set}.
Therefore, let $x>0$, then it follows by Definition~\ref{def:cusp} that 
\begin{align*}
x=\int_0^x dy & =-\int_0^x \frac{\sqrt{s-\phi(y)}}{\sqrt{M-\phi(y)}(\phi(y)-m)} 
\phi_x(y)dy\\
&= \int_{\phi(x)}^s\frac{\sqrt{s-z}}{\sqrt{M-z}(z-m)}dz\\ 
& \geq \frac{1}{(M-m)^{3/2}} \int_{\phi(x)}^s \sqrt{s-z}dz\\
& = \frac{2}{3}\left(\frac{s-\phi(x)}{M-m}\right)^{3/2}
\end{align*}
or equivalently
\[ 
0\leq s-\phi(x)=\phi(0)-\phi(x)\leq (M-m) \left(\frac32 x\right)^{2/3} \quad 
\text{ for all } x\in\Real
\]
since $\phi(x)$ is symmetric. Hence, given some $\varepsilon >0$, there exists a 
$\delta>0$ such that 
\[
\vert \phi(0)-\phi(x)\vert <\varepsilon \quad \text{ for all } \vert x\vert 
<\delta.
\]
Moreover, it is well-known that $y_0(\xi)$ is Lipschitz continuous with 
Lipschitz constant at most one. Thus 
\[ 
\vert y_0(\xi)-y_0(\bar\xi)\vert <\delta \quad \text{ for all } \vert 
\xi-\bar\xi\vert <\delta
\]
and the definition of $y_0(\xi)$ yields
\begin{align*}
\vert \bar\xi-\xi\vert & =\vert y_0(\bar \xi)-y_0(\xi)\vert +\vert 
\int_{y_0(\xi)}^{y_0(\bar \xi)} \phi_x^2(z)dz\vert \\
& \geq \vert y_0(\bar\xi)-y_0(\xi)\vert +\vert \int_{y_0(\xi)}^{y_0(\bar\xi)} 
\frac{(M-s)(s-\varepsilon-m)^2}{\varepsilon} dz\vert\\
& \geq \vert y_0(\bar \xi)-y_0(\xi)\vert 
\left(\frac{\varepsilon+(M-s)(s-\varepsilon-m)^2}{\varepsilon}\right)
\end{align*}
Hence 
\[
\frac{y_0(\bar \xi)-y_0(\xi)}{\bar\xi-\xi}\leq 
\frac{\varepsilon}{\varepsilon+(M-s)(s-\varepsilon-m)^2} \quad
\text{ for } \vert \bar\xi-\xi\vert <\delta.
\]
Since we can choose $\varepsilon$ to be any positive real number, we have 
\[ 
y_{0,\xi}(\bar\xi)=\lim_{\xi\to\bar\xi} \frac{y_0(\bar 
\xi)-y_0(\xi)}{\bar\xi-\xi}=0.
\]

{\bf The cusp is not traveling along a single characteristic.}
Finally, we want to show that the cusp is not travelling along a characteristic 
in Lagrangian coordinates. In particular, we are going to show that the peak, 
metaphorically speaking, jumps from one characteristic to the next. Recall 
therefore that the cuspon with exponential decay and speed $s+\kappa$ is given 
by 
\[ 
u(t,x)=\phi(x-(s+\kappa)t)+\kappa.
\]
Thus by definition we have 
\begin{align*}
P(t,x)&=\frac14 \int_\Real e^{-\vert x-z\vert} (2u^2+u_x^2)(t,z)dz\\
& = \frac14 \int_\Real e^{-\vert (x-(s+\kappa)t)-(z-(s+\kappa)t)\vert} 
(2(\phi(z-(s+\kappa)t)+\kappa)^2+\phi_x^2(z-(s+\kappa)t))dz\\
&=\frac14 \int_\Real e^{-\vert 
(x-(s+\kappa)t)-z\vert}(2(\phi(z)+\kappa)^2+\phi_x^2(z))dz\\
& = P(0,x-(s+\kappa)t).
\end{align*}
Similar considerations yield
\[
P_x(t,x)=P_x(0,x-(s+\kappa)t).
\]
Moreover, $P(0,x)$ is an even function, and hence $P_x(0,x)$ is an odd function 
and satisfies
\[
P_x(t,(s+\kappa)t)=P_x(0,0)=0.
\]
Denote by $y(t,\xi)$ the characteristics at time $t$, where $y(0,\xi)$ is 
defined through \eqref{eq:Ldef1} and satisfies
\[
y_t(t,\xi)=U(t,\xi).
\]
Then we have 
\begin{subequations}
\begin{align}\label{travel:func}
U(t,\xi)&=u(t,y(t,\xi))=\phi(y(t,\xi)-(s+\kappa)t)+\kappa,\\
P(t,\xi)& =P(t,y(t,\xi))=P(0,y(t,\xi)-(s+\kappa)t)\\
Q(t,\xi)&= P_x(t,y(t,\xi))=P_x(0,y(t,\xi)-(s+\kappa)t)
\end{align}
\end{subequations}
and hence 
\begin{equation}
(y(t,\xi)-(s+\kappa)t)_t=\phi(y(t,\xi)-(s+\kappa)t)-s\leq 0.
\end{equation}
Thus the characteristics to the left and to the right of the cusp travel at a 
speed slower than the one of the cusp. In particular, one can show using 
subsolutions, that to any $\xi>\bar\xi$ there exists a time $T>0$ such that 
$y(t,\xi)>(s+\kappa)t$ for all $t<T$, $y(T,\xi)=(s+\kappa)T$ and 
$U(t,\xi)=(s+\kappa)$. Obviously we have for $\xi>\bar\xi$ that 
$y(t,\xi)<(s+\kappa)t$ for all $t>0$.

Let $\tilde\xi\in\Real$ such that $y(T,\tilde\xi)=(s+\kappa)T$, then we aim at showing that $y(t,\tilde\xi)<(s+\kappa)t$ for all $t>T$. Or equivalently, we can show for $z(t,\tilde\xi)=y(t,\tilde\xi)-(s+\kappa)t$, that $z(t,\tilde\xi)<0$ for all $t>T$. Therefore note that 
\begin{align}
 z_t(T,\tilde\xi)& = U(T,\tilde\xi)-(s+\kappa)=\phi(z(T,\tilde\xi))-s=0,\\
 z_{tt}(T,\tilde\xi) & = -Q(T, \tilde\xi)=-P_x(0,z(T,\tilde\xi))=0\\
 z_{ttt}(T,\tilde\xi)& = -Q_t(T,\tilde\xi)=?
\end{align}
Hence, if we can show that $Q_t(T,\tilde\xi)$ exists and is positive, it follows that the function $U(t,\tilde\xi)$, seen as a function of $t$ attains a maximum at $t=T$. Thus $z_t(t,\tilde\xi)<0$ for all $t>T$ and, in particular, $z(t,\tilde\xi)<0$ for all $t>T$.

Recall that $z(t,\xi)$ is bijective and continuous and that $Q(t,\xi)=P_x(0, z(t,\xi))$. Thus for all $\xi\in\Real\backslash\{\tilde\xi\}$, we can apply the chain rule to $Q(t,\xi)$ and obtain 
\begin{align}
 Q_t(T,\xi)& = P_{xx}(0,z(T,\xi))z_t(T,\xi)\\
 & = (P(0,z(T,\xi))-(\phi(z(T,\xi))-\kappa)^2)(\phi(z(T,\xi))-s)\\
 & \quad \phantom{space}+(M-\phi(z(T,\xi)))(\phi(z(T,\xi))-m)^2.
\end{align}
 If $Q_t(T,\xi)$ is continuous as a function of $\xi$, we get that 
 \begin{equation}
 Q_t(T,\tilde\xi)=(M-s)(s-m)^2>0, 
 \end{equation}
since $P(0,x)$ and $\phi(x)$ are uniformly bounded and $z_t(T,\tilde\xi)=0$. 

Thus establishing the existence and continuity of $Q_t(t,\xi)$ will finish the proof of the claim. 
By definition,
\begin{align}
 Q(t,\xi)&=-\frac14 \int_{-\infty}^\xi e^{-(\xi-\eta)}e^{-(\zeta(t,\xi)-\zeta(t,\eta))}(2U^2y_\xi(t,\eta)+h(t,\eta))d\eta\\
 &\quad + \frac14 \int_\xi^\infty e^{-(\eta-\xi)}e^{-(\zeta(t,\eta)-\zeta(t,\xi))}(2U^2y_\xi(t,\eta)+h(t,\eta))d\eta,
\end{align}
where we used $y(t,\xi)=\xi+\zeta(t,\xi)$.
Since the functions 
\begin{subequations}
 \begin{align}
  f_1(t,\xi,\eta)& = e^{-(\zeta(t,\xi)-\zeta(t,\eta))}(2U^2(t,\eta)y_\xi(t,\eta)+h(t,\eta)) \\
  f_2(t,\xi,\eta)&= e^{-(\zeta(t,\eta)-\zeta(t,\xi))}(2U^2(t,\eta)y_\xi(t,\eta)+h(t,\eta))
 \end{align}
\end{subequations}
are both differentiable with respect to $t$ and $f_{1,t}(t,\xi,\eta)$ and $f_{2,t}(t,\xi,\eta)$ can be uniformly bounded for $\xi$, $\eta\in\Real$ and $t$ on a finite time interval, \cite[Theorem 2.27]{Folland} implies the existence of $Q_t(t,\xi)$ and that 
\begin{align*}
 Q_t(t,\xi)& =\frac14\int_\Real e^{-\vert y(t,\xi)-y(t,\eta)\vert}(U(t,\xi)-U(t,\eta))(2U^2y_\xi(t,\eta)+h(t,\eta)) d\eta\\
 & \quad - \frac14 \int_\Real \sign(\xi-\eta) e^{-\vert y(t,\xi)-y(t,\eta)\vert }
\\
 & \qquad \phantom{space}\times(-4QU(t,\eta) y_\xi(t,\eta)+4U^2U_\xi(t,\eta)-2PU_\xi(t,\eta))d\eta.
\end{align*}
Finally, following the same lines as the proof of \cite[Lemma 3.1]{GHR2}, the continuity of $Q_t(t,\xi)$ with respect to $\xi$ can be established.

\vspace{0.5cm}

\section{Observation}\label{Observation}

We had now a closer look at two particular initial data for the CH equation, 
which yield on the one hand a solution, with accumulating breaking times, see 
Section~\ref{sec:ex} and on the other hand a solution which enjoys wave breaking 
at any time, but the associated measure has neither a discrete nor a singular 
continuous part, see Section~\ref{sec:cusp}. Thus naturally the question arises, 
can we find to some initial data $(u_0,\mu_0)\in\D$ a global conservative 
solution $(u(t,x),\mu(t,x))\in\D$ such that there exist $0<t_1<t_2$ such that 
$\supp(\mu_d(t,x))\not=\emptyset$ for all $t\in[t_1,t_2]$?

Reformulating this question in Lagrangian coordinates yields. Can we find some 
initial data $X_0=(y_0,U_0,h_0)\in\F$ with solution $X(t)=(y(t),U(t),h(t))$ such 
that there exists $0<t_1<t_2$ such that $y_\xi(t,\xi)=0$ on an interval of 
positive length for all $t\in[t_1,t_2]$? 

Due to \cite[Lemma 2.7 (ii)]{HR}, which we state in a moment for the sake of 
completeness, the answer is no.

\begin{lemma}[{\cite[Lemma 2.7 (ii)]{HR}}]
Given initial data $X_0=(\zeta_0, U_0, h_0)$ in $\mathcal{F}$, let 
$X(t)=(\zeta(t),U(t),h(t))$ be the short-time solution of \eqref{eq:chsyseq} in 
$C([0,T],\mathcal{F})$ for some $T>0$ with initial data $X_0=(\zeta_0, U_0, 
h_0)$. Then for almost every $t\in [0,T]$, $y_\xi(t,\xi)>0$ for almost every 
$\xi\in\Real$.
\end{lemma}

\end{document}